\documentclass[a4paper,leqno,11pt]{article}
\usepackage{amsmath, amssymb, amsthm, amsfonts}
\usepackage{mathrsfs}

\usepackage{color,psfrag}

\theoremstyle{plain}
\newtheorem{theorem}{Theorem}

\newtheorem{lemma}{Lemma}

\theoremstyle{remark}
\newtheorem{remark}{Remark}

\renewcommand{\Re}{{\rm Re\,}}
\renewcommand{\Im}{{\rm Im\,}}

\numberwithin{equation}{section}

\title{Average of Hardy's function at Gram points}
\author{Xiaodong Cao, Yoshio Tanigawa and Wenguang Zhai\thanks{Wenguang Zhai is supported by 
National Natural Science Foundation of China (Grant No. 11971476).}}
\date{} %\today}

\begin{document}

\maketitle

\begin{abstract}
Let $Z(t)=\chi^{-1/2}(1/2+it)\zeta(1/2+it)=e^{i\theta(t)}\zeta(1/2+it)$ be Hardy's function and
$g(n)$ be the $n$-th Gram points defined by $\theta(g(n))=\pi n$.  Titchmarsh proved that $\sum_{n \leq N} Z(g(2n))
=2N+O(N^{3/4}\log^{3/4}N) $ and $\sum_{n \leq N} Z(g(2n+1)) =-2N+O(N^{3/4}\log^{3/4}N)$.
We shall improve the error terms to $O(N^{1/4}\log^{3/4}N \log\log N)$.

\end{abstract}

%***********************************************************

\section{Introduction}

%***********************************************************

Let $s=\sigma+it$ be a complex variable and let $\zeta(s)$ be the Riemann zeta-function. 
It satisfies the functional equation
\begin{equation} \label{FE}
\zeta(s)=\chi(s)\zeta(1-s),
\end{equation}
where
\begin{equation} \label{def-chi}
\chi(s)=2^s\pi^{s-1}\sin\left(\frac{\pi s}{2}\right)\Gamma(1-s)=\pi^{s-1/2}\frac{\Gamma(\frac12-\frac{s}{2})}{\Gamma(\frac{s}{2})}
\end{equation}
and the asymptotic behavior of $\chi(s)$ is given by
\begin{equation} \label{chi-asymp}
\chi(s)=\left(\frac{|t|}{2\pi}\right)^{1/2-\sigma-it}e^{i(t \pm \frac{\pi}{4})}
 \left(1+O\left(\frac{1}{|t|}\right)\right),  \quad (|t| \geq 1)
\end{equation}
with $t \pm \frac{\pi}{4}=t+{\rm sgn}(t)\frac{\pi}{4}$  (see Ivi\'c \cite[(1.25)]{I0}).
Hardy's function $Z(t)$ is defined by
\begin{equation} \label{def_Hardy}
Z(t)= \chi^{-1/2}\left(\frac12+it\right)\zeta\left(\frac12+it\right),
\end{equation}
From \eqref{FE}, \eqref{def-chi} and \eqref{def_Hardy}, it follows that $Z(t)$ is a real-valued even function
for real $t$ and $|Z(t)|=|\zeta(1/2+it)|$. Thus the real zeros of $Z(t)$ coincide with the zeros of
$\zeta(s)$ on the critical line $\Re s=1/2$. 
Furthermore we have an equivalent expression of \eqref{def_Hardy} 
$$
Z(t)=e^{i\theta(t)}\zeta\left(\frac12+it\right),
$$
where
\begin{align}  \label{def-theta}
\theta(t)&=-\frac{1}{2i} \log\chi\left(\frac12+it\right) \\
&=\Im\left(\log \Gamma\left(\frac14+\frac{it}{2}\right)\right)-\frac{t}{2}\log \pi
\in \mathbb{R}.  \nonumber 
\end{align}
It is well-known that the function $\theta(t)$ is strictly monotonic increasing for $t \geq 6.5$.  (cf. \eqref{theta1-seikaku} below.)

Historically, Hardy proved the infinity of numbers of zeros
of $\zeta(s)$ on the critical line in 1914. A little later Hardy and Littlewood gave another proof by 
%comparing the integral of $Z(t)$ and $|Z(t)|$.
showing that $\int_0^T Z(t)dt \ll T^{7/8}$ and $\int_0^T|Z(t)|dt \gg T$  (see Chandrasekharan \cite[Chapter II, $\S$ 4 and Notes 
on Chapter II]{C} 
or Titchmarsh \cite[10.5]{T}).  
It must be mentioned that the mean value estimate of $Z(t)$ was improved extensively to $\int_0^TZ(t) \ll
T^{1/4+\varepsilon}$ by Ivi\'c \cite{I2} in 2004, where $\varepsilon$ is an arbitrary small positive number.
%On his result it can be said that the interest on Hardy's function has revived. 
See Ivi\'c's monograph \cite{I3} for the recent development of the theory of Hardy's function.

Before Hardy, Gram calculated zeros of $\zeta(1/2+it)$ and observed that
the points $t$ such that $\Re \zeta(1/2+it) \in \mathbb{R}$ and the zeros of $\zeta(1/2+it)$ are distributed
alternately.
There is also tendency that $\Re \zeta(1/2+it)$ takes positive values and $\Im \zeta(1/2+it)$ takes positive
and negative values regularly.  See e.g. the graphs of $\zeta(1/2+it)$ in Akiyama and Tanigawa \cite{AT}.

%From \eqref{FE} and \eqref{def-chi} we have
%$$
%Z(t)=e^{i\theta(t)}\zeta\left(\frac12+it\right),
%$$
%where
%\begin{equation}  \label{def-theta}
%\theta(t)=-\frac{1}{2i} \log\chi\left(\frac12+it\right)=\Im\left(\log \Gamma\left(\frac14+\frac{it}{2}\right)\right)-\frac{t}{2}\log \pi
%\in \mathbb{R}.
%\end{equation}
%It is well-known that the function $\theta(t)$ is strictly monotonic increasing for $t \geq 6.5$.  (cf. \eqref{theta1-seikaku} below.)

For $n \geq -1$, let $g(n)>7 $ be the $n$-th Gram point defined by
$$
\theta(g(n))=\pi n.
$$
Obviously
$$
\zeta\left(\frac12+ig(n)\right)=(-1)^n Z(g(n)).
$$
Gram's law is stated that there exists a zero of $Z(t)$ for some $t\in [g(n), g(n+1)]$.
The first twelve Gram points are (Haselgrove and Miller \cite{HM})
\begin{align*}
&g(-1)=9.6 \ldots, && g(0)=17.8\ldots, && g(1)=23.1\ldots, && g(2)=27.6\ldots, \\
&g(3)= 31.7\ldots, && g(4)=35.4 \ldots, && g(5)=38.9\ldots, && g(6)=42.3\ldots, \\
&g(7)=45.5\ldots, && g(8)=48.7\ldots, && g(9)=51.7\ldots, && g(10)=54.7\ldots.
\end{align*}
At present it is known that there is a positive proportion of failures of Gram's law (see e.g. Trudgain \cite{Tru}
and Ivi\'c \cite[p.~112]{I3}).

As for the distribution of $Z(g(n))$ on the average, Titchmarsh showed that, for a fixed large integer $ M $,
\begin{align}
\sum_{M+1}^{N} Z(g(2n))&=2(N-M)+O(N^{3/4}\log^{3/4}N), \label{Tit-1} \\
\sum_{M+1}^{N} Z(g(2n+1))&=-2(N-M)+O(N^{3/4}\log^{3/4}N)  \label{Tit-2}
\end{align}
\cite[10.6]{T}\footnote{Titchmarsh uses the notation $t_\nu$ for the Gram points $g(\nu)$.  
In \cite[Theorem 10.6, p.~263]{T} he asserted that $\sum_{\nu=\nu_0}^{N}Z(t_{2\nu}) \sim 2N$ and 
$\sum_{\nu=\nu_0}^{N}Z(t_{2\nu+1}) \sim -2N$, but in fact he obtained the error terms as in \eqref{Tit-1} and \eqref{Tit-2}, see p.264 of \cite{T}.},
where he used the approximation 
\begin{equation} \label{approx}
Z(g(n))=2(-1)^n \sum_{m \leq \sqrt{g(n)/2\pi}} m^{-1/2} \cos(g(n)\log m)+O(g(n)^{-1/4})
\end{equation}
obtained by the classical approximate functional equation of the Riemann zeta-function due to Hardy and Littlewood.
In Ivi\'c \cite[Theorem 6.5]{I3}, the error terms of \eqref{Tit-1} and \eqref{Tit-2} are improved to $O(N^{3/4}\log^{1/4}N)$.
In \cite{T0} Titchmarsh also proved
\begin{align}  \label{T-IV}
& \sum_{\nu=M+1}^{N}\left(\sum_{n=1}^{[\sqrt{g(\nu)/2\pi}]}\frac{\cos(g(\nu)\log n)}{\sqrt{n}}\right)
=N+O(N^{1/4} \log^{-1/4}N), \\
& \sum_{n \leq N}Z(g(n))Z(g(n+1))=-2(\gamma+1)N+o(N),  \nonumber
\end{align}
where $\gamma$ is Euler's constant and conjectured
$$
\sum_{n \leq N}Z(g(n))^2 Z(g(n+1))^2 \ll N \log^A N
$$
with some positive constant $A$.  This conjecture was proved by Moser in \cite{M} (See Ivi\'c \cite[Notes of Chapter 6]{I3}).

The purpose of the present paper is to show the following theorem.
\begin{theorem} \label{thm1}
 We have
\begin{align}
\sum_{n \leq N}Z(g(2n))&=2N+O\left(N^{1/4}\log^{3/4}N \log\log N \right), \label{evencase}  \\
\sum_{n \leq N}Z(g(2n+1))&=-2N+O\left(N^{1/4}\log^{3/4}N \log \log N \right). \label{oddcase}
\end{align}
\end{theorem}

\begin{remark}
%In \cite{T0}, Titchmarsh proved that ($t_{\nu}=g(\nu)$ in our notation)
%\begin{align} \label{T-IV}
%\sum_{\nu=M+1}^{N}\left(\sum_{n=1}^{[\sqrt{t_{\nu}/2\pi}]}\frac{\cos(t_{\nu}\log n)}{\sqrt{n}}\right)
%=N+O(N^{1/4} \log^{-1/4}N).
%\end{align}
Note that there is not sign $(-1)^n$ in the sum on the left hand side of \eqref{T-IV}.
It seems that our Theorem 1 does not follow from \eqref{T-IV}, since when we consider $g(2n)$ and $g(2n+1)$ separately
we may not be able to use the first derivative test directly (cf. Trudgian \cite{Tru}).
\end{remark}

%**********************************************

\section{The function $\theta(t)$}

%**********************************************

As in (1.24)--(1.26) in Ivi\'c \cite{I3}, the function $\theta(t)$ defined by \eqref{def-theta} and its derivatives have
asymptotic expansions. In particular
\begin{align}
\theta(t)&=\frac{t}{2}\log\frac{t}{2\pi}-\frac{t}{2}-\frac{\pi}{8}+\frac{1}{48t}+\frac{7}{5760t^3}+O(t^{-5}),\label{theta}\\
\theta'(t)&=\frac12\log\frac{t}{2\pi}+O(t^{-2}),  \label{theta1-asymp} \\
\theta''(t)&=\frac{1}{2t}+O(t^{-3}),   \label{theta2-asymp}
\end{align}
(\cite[(1.27)]{I3}) and
$$
\theta'''(t) \ll \frac{1}{t^2},  \quad  \theta^{(4)}(t) \ll \frac{1}{t^3}.
$$
For $\theta'(t)$ and $\theta''(t)$ we shall need more precise formulas in the proof Theorem \ref{thm1}.

\begin{lemma} \label{lem-4}
For $t \geq 6$ we have
\begin{align}
\theta'(t)&=\frac12\log\frac{t}{2\pi}-\frac{1}{48t^2}+V_1(t) \quad \text{with \  $|V_1(t)| \leq 0.07 \, t^{-3}$} \label{theta1-seikaku} \\
\intertext{and}
\theta''(t)&=\frac{1}{2t}+V_2(t) \quad \text{with \ $|V_2(t)| \leq 0.46 \, t^{-3}$}.  \label{theta2-seikaku}
\end{align}
\end{lemma}

\proof
Let $B_{n}(x)$ be the $n$-th Bernoulli polynomial defined by $\frac{te^{xt}}{e^t-1}=\sum_{n=0}^{\infty}\frac{B_{n}(x)}{n!}t^n$,
$(|t|<2\pi)$ and $B_n=B_n(0)$ be the $n$-th Bernoulli number.
It is well known that for $|\arg z|<[\pi|$,
\begin{align*}
\log \Gamma(z)&=\left(z-\frac12\right)\log z-z+\frac12 \log 2\pi + \sum_{r=1}^n \frac{B_{2r}}{2r(2r-1)}z^{-2r+1}\\
& \quad -(2n)!\int_0^{\infty}\frac{P_{2n+1}(x)}{(x+z)^{2n+1}}dx,
\end{align*}
where $P_n(x)$ is defined by $P_n(x)=\frac{B_n(x-[x])}{n!}$, see Wang and Guo \cite[p.114 (8)]{WG}.
Differentiating the above formula we have
\begin{align}  \label{digamma}
\psi(z)=\log z-\frac{1}{2z}-\sum_{r=1}^n \frac{B_{2r}}{2r}z^{-2r}+(2n+1)!\int_0^{\infty}\frac{P_{2n+1}(x)}{(x+z)^{2n+2}}dx.
\end{align}
Here $\psi(z)$ denotes the digamma function: $\psi(z)=\frac{\Gamma'(z)}{\Gamma(z)}$.
Let $z=\sigma+it$, ($0<\sigma<1$ and $t>0$). Then using $|P_{m}(x)| \leq \frac{4}{(2\pi)^m}$ for $m \geq 1$ (\cite[p. 11 (15)]{WG},
we have
\begin{align*}
\left|\int_0^{\infty}\frac{P_{2n+1}(x)}{(x+z)^{2n+2}}dx\right| & \leq \frac{4}{(2\pi)^{2n+1}}
\int_0^\infty \frac{1}{((x+\sigma)^2+t^2)^{n+1}}dx \\
& \leq \frac{4}{(2\pi)^{2n+1}}t^{-2n-1} \int_0^{\infty}\frac{1}{(u^2+1)^{n+1}}du\\
&=\frac{4}{(2\pi)^{2n+1}}\frac{(2n-1)!!}{(2n)!!}t^{-2n-1}.
\end{align*}
Take $n=1$ in \eqref{digamma}, then for $z=\sigma+it, \ (0<\sigma<1, t>0)$ we have 
$$
\psi(z)=\log z-\frac{1}{2z}-\frac{1}{12z^2}+K_1, \quad |K_1| \leq \frac{3}{4\pi^2}t^{-3}.
$$
Now we differentiate the both side of \eqref{def-theta} and get
\begin{align*}
\theta'(t)&=\frac14\left(\psi\left(\frac14+\frac{it}{2}\right)+\psi\left(\frac14-\frac{it}{2}\right)\right)-\frac12\log \pi \\
&=\frac14\left(\log\left(\frac{1}{16}+\frac{t^2}{4}\right)-\frac14\frac{1}{\frac{1}{16}+\frac{t^2}{4}}
-\frac{1}{24}\frac{-t^2+\frac14}{(\frac{1}{16}+\frac{t^2}{4})^2}+K_1+\bar{K}_1\right)-\frac12\log\pi.
\end{align*}
Since
%\begin{equation}
\begin{gather*}
-\frac{1}{32t^4} < \log\left(\frac{1}{16}+\frac{t^2}{4}\right)-2\log\frac{t}{2}-\frac{1}{4t^2}<0, \\
-\frac{1}{t^4}< \frac{1}{\frac{1}{16}+\frac{t^2}{4}}-\frac{4}{t^2}<0, \quad
0< \frac{-t^2+\frac14}{(\frac{1}{16}+\frac{t^2}{4})^2}+\frac{16}{t^2}<\frac{12.4}{t^4},
\end{gather*}
%\end{equation}
we have
\begin{equation*}
\theta'(t)=\frac12\log\frac{t}{2\pi}-\frac{1}{48t^2}+V_1  \quad  \text{with \ $-0.061 \, t^{-3}< V_1 < 0.0485 \, t^{-3}$}.
\end{equation*}
This proves the assertion \eqref{theta1-seikaku}.

The assertion \eqref{theta2-seikaku} is proved in a similar way.

\endproof

In connection with the definition of the Gram point $g(n)$,  we define a function $g(x)$ for real variable $x \geq -1$ by
\begin{equation}  \label{gramfunction}
\theta(g(x))=\pi x.
\end{equation}
The function $g(x)$ is uniquely determined and monotonic increasing.
When $x$ is an integer, $g(x)$ coincides with the definition of Gram points.
For the order of $g(n)$ it is known that
\begin{align*}
&g(n)=\frac{2\pi n}{\log n}\left\{1+ \frac{1+\log \log n}{\log n}+O\left(\left(\frac{\log\log n}{\log n}\right)^2\right)\right\}  \\
&n=\frac{g(n)}{2\pi} \log g(n) \left\{1-\frac{\log 2\pi e}{\log g(n)}-\frac{\pi}{4g(n)\log g(n)}+O\left(\frac{1}{g(n)^2\log g(n)}\right)\right\}.
\end{align*}
See Ivi\'c \cite[Theorem 6.1]{I3}, also Bruijin \cite{B}. Note that they hold for any positive numbers $n$.

From this definition of $g(x)$, we have
\begin{align*}
g'(x)&=\frac{\pi}{\theta'(g(x))},\\
g''(x)&=-\frac{\pi^2 \theta''(g(x))}{\theta'(g(x))^3}, \\
g'''(x)&=-\pi^3 \frac{\theta'''(g(x))\theta'(g(x))-3\theta''(g(x))^2}{\theta'(g(x))^5}, \\
g^{(4)}(x)&=-\pi^4 \frac{\theta^{(4)}(g(x))\theta'(g(x))^2-10\theta'''(g(x))\theta''(g(x))\theta'(g(x))+15\theta''(g(x))^3}{\theta'(g(x))^7}.
\end{align*}

\vskip 1cm

%***********************************************************

\section{Some Lemmas}

%*************************************************************
%For \cite{T0} and \cite{I3}, the approximation
%$$
%Z(t)=2 \sum_{n \leq \sqrt{t/2\pi}}n^{-1/2}\cos(\theta(t)-t\log n)+O(t^{-1/4})
%$$
% is used. However,
The error term of \eqref{approx} is rather big for our purpose. So we apply an approximation of $Z(t)$ with smooth weight
which is due to Ivi\'c \cite{I3}.
\begin{lemma} \label{lem-1}
We have
\begin{equation} \label{AFE}
Z(t)=2\sum_{m=1}^{\infty} \frac{1}{m^{\frac12}}\rho\left(m\sqrt{\frac{2\pi}{t}}\right)\cos(\theta(t)-t\log m)+O(t^{-5/6}),
\end{equation}
where $\rho(t)$ is a real-valued function such that
\begin{enumerate}
\item[{\rm (i)}]  $\rho(t) \in C^{\infty}(0,\infty)$,
\item[{\rm (ii)}]  $\rho(t)+\rho(1/t)=1$ for $t>0$,
\item[{\rm (iii)}]  $\rho(t)=0$ for $t \geq 2$.
\end{enumerate}
\end{lemma}

\proof 
This lemma is obtained from the definition of $Z(t)$ and Theorem 4.16 of Ivi\'c \cite{I3}. See also (4.81) of \cite{I3}.
\endproof

In \cite[Lemma 4.15]{I3}, Ivi\'c constructed a function $\rho(t)$ in more general form in such a way that, instead of (iii),  it 
satisfies $\rho(t)=0$ for $t \geq b$ for any fixed $b>1$.
But the choice $b=2$ is sufficient for our purpose. More explicitly it is given as follows with the choices $\alpha=3/2,
\beta=1/2$ in \cite[Lemma 4.15]{I3}. Let
$$
\varphi(t)=\exp((t^2-1/4)^{-1}) \left\{\int_{-\frac12}^{\frac12}\exp((u^2-1/4)^{-1})du \right\}^{-1}
$$
if $|t| <\frac12$ and $\varphi(t)=0$ if $|t| \geq \frac12$, and let
$$
f(x)=\int_{x-\frac32}^{x+\frac32}\varphi(t)dt=\int_{-\infty}^x\left(\varphi(t+3/2)-\varphi(t-3/2)\right)dt.
$$
The function $f(x)$ is infinitely differentiable in $(-\infty, \infty)$, $f(x) \geq 0$ and
$$
f(x)= \begin{cases} 0  &  \text{for $|x| \geq 2$} \\
                    1  &  \text{for $|x| \leq 1$}
       \end{cases}
$$
(see \cite[p.88]{I3} for details).  Then the function
\begin{equation} \label{Ivic-rho}
\rho(x)=\frac12\left(1+f(x)-f\left(\frac{1}{x}\right)\right)
\end{equation}
satisfies (i), (ii) and (iii) of Lemma \ref{lem-1}.

\begin{lemma} \label{lem-2}
Let $\rho(t)$ be the function as above. When $x \to 0$, we have
\begin{equation} \label{rhokyodo}
\rho(1+x)=\frac12+O(|x|^C),
\end{equation}
where $C$ is any positive large constant.
\end{lemma}

\proof  This property depends on the construction as above.
%We shall use the explicit construction of $\rho(t)$ in Lemma \ref{lem-1} due to Ivi\'c.
%According to Lemma 4.15 of Ivi\'c \cite{I3} with the choices $\alpha=3/2, \beta=1/2$, let
%$$
%\varphi(t)=\exp((t^2-1/4)^{-1}) \left\{\int_{-\frac12}^{\frac12}\exp((u^2-1/4)^{-1})du \right\}^{-1}
%$$
%if $|t| <\frac12$ and $\varphi(t)=0$ if $|t| \geq \frac12$, and let
%$$
%f(x)=\int_{x-\frac32}^{x+\frac32}\varphi(t)dt=\int_{-\infty}^x\left(\varphi(t+3/2)-\varphi(t-3/2)\right)dt.
%$$
%Note that $f(x) \in C^{\infty}(-\infty, \infty)$, $f(x) \geq 0$ and
%$$
%f(x)= \begin{cases} 0  &  \text{for $|x| \geq 2$} \\
%                    1  &  \text{for $|x| \leq 1$}
%       \end{cases}
%$$
%(see \cite[p.88]{I3} for details).
%Then the function
%\begin{equation} \label{Ivic-rho}
%\rho(x)=\frac12\left(1+f(x)-f\left(\frac{1}{x}\right)\right)
%\end{equation}
%satisfies (i), (ii) and (iii) of Lemma \ref{lem-1}.

When $x$ is positive and small, we have
$$
f(1+x)=\int_{-\frac12+x}^{\frac12}\varphi(t)dt.
$$
Using $\int_{-\frac12}^{\frac12}\varphi(t)dt=1$ we find that
\begin{align*}
1-f(1+x) =\int_{-\frac12}^{-\frac12+x} \varphi(t)dt \leq x \max \left\{\varphi(t) \  \Big | \  -\frac12 \leq t \leq -\frac12+x\right\}
\ll x^{C},
\end{align*}
where $C$ is any positive large constant. On the other hand $f(\frac{1}{1+x})=1$.
Thus from \eqref{Ivic-rho} it follows that
$$
\rho(1+x)=\frac12+O(x^C)
$$
when $ x \to 0+$. 

It is proved by the same way when $x \to 0-$.

\endproof

We need the following two lemmas on exponential sums.

\begin{lemma} \label{lem-fdt}
Let $f(x)$ and $\varphi(x)$ be real-valued functions which satisfy the following conditions on the interval $[a,b]$.

$(1)$ \  $f'(x)$ is continuous and monotonic on $[a,b]$ and $|f'(x)| \leq \delta <1$.

$(2)$ \  $\varphi(x)$ is positive monotonic and $\varphi'(x)$ is continuous,  and 
there exist numbers $0<H, \ 0<b-a \leq  U$ such that 
$$
\varphi(x) \ll H, \quad \varphi'(x) \ll HU^{-1}.
$$ 
Then we have
\begin{equation} \label{exp-sum-1}
\sum_{a<n \leq b} \varphi(n)\exp(2\pi i f(n))=\int_a^b \varphi(x)\exp(2\pi i f(x))dx+O\left(\frac{H}{1-\delta}\right).
\end{equation}
Furthermore, if $0<\delta_1<f'(x) \leq \delta<1 $ on $[a,b]$, then we have
\begin{equation} \label{exp-sum-2}
\sum_{a<n \leq b} \varphi(n)\exp(2\pi i f(n)) \ll \frac{H}{\delta_1}+\frac{H}{1-\delta}.
\end{equation}
\end{lemma}

\proof The assertion \eqref{exp-sum-1} is obtained by Lemma 1.2 of Ivi\'c \cite{I0} and partial summation.
The second assertion \eqref{exp-sum-2} is obtained by applying (2.3) of Ivi\'c \cite{I0} 
(the so-called first derivative test) on the integral of \eqref{exp-sum-1}.  
See also Karatsuba and Voronin \cite[p.~70 Corollary 1]{KV}. 
\endproof

\begin{lemma}[{Karatsuba and Voronin \cite[Chapter III, Theorem 1]{KV}}] \label{lem-3}
Suppose that the real-valued function $\varphi(x)$ and $f(x)$ satisfy the following conditions on the interval $[a,b]$:
\begin{enumerate}
\item[$(1)$] $f^{(4)}(x)$ and $\varphi''(x)$ are continuous;
\item[$(2)$] there exist numbers $H,\, U,\, A, \, 0<H, 1 \ll A \ll U, \, 0<b-a \leq U$, such that
\begin{align*}
& A^{-1} \ll f''(x) \ll A^{-1},  && f^{(3)}(x) \ll A^{-1}U^{-1}, && f^{(4)}(x) \ll A^{-1} U^{-2} \\
& \varphi(x) \ll H, && \varphi'(x) \ll HU^{-1}, && \varphi''(x) \ll HU^{-2}.
\end{align*}
\end{enumerate}
Suppose that the numbers $x_n$ are determined from the  equation
$$
f'(x_n)=n.
$$
Then we have
\begin{equation}  \label{exp-sum-trans}
\sum_{a<x \leq b} \varphi(x)\exp(2\pi i f(x))=\sum_{f'(a) \leq n \leq f'(b)}c(n)W(n)+R,
\end{equation}
where
\begin{align}
&R=O\left(H(A(b-a)^{-1}+T_a+T_b+\log(f'(b)-f'(a)+2))\right); \label{gosaR} \\[1ex]
&T_{\mu}=\left\{ \begin{array}{cl} 0, & \text{if $f'(\mu)$ is an integer}, \\
                                 \min\left(\| f'(\mu) \| ^{-1}, \sqrt{A}\right), & \text{if $\| f'(\mu) \| \neq 0$};
                \end{array} \right.  \label{def-T} \\[1ex]
&c(n)=\left\{\begin{array}{cl} 1, & \text{if $f'(a)<n<f'(b)$,} \\
                              1/2, & \text{if $n=f'(a)$ or $n=f'(b)$};
              \end{array} \right. \nonumber \\[1ex]
&W(n)=\frac{1+i}{\sqrt{2}}\frac{\varphi(x_n)}{\sqrt{f''(x_n)}}\exp(2\pi i(f(x_n)-nx_n)).  \label{shukou}
\end{align}
\end{lemma}

\bigskip

\begin{remark} \label{remark-KV}
For the assertion of Lemma \ref{lem-3}, the condition $A \ll U$ in (2) is not necessary.
See also Jia \cite[Lemma 5]{J}.
\end{remark}

\begin{remark}  \label{remark-KV2}
When $f''(x)$ is negative and $A^{-1} \ll -f''(x) \ll A^{-1}$, the sum on the right hand side of \eqref{exp-sum-trans} should be replaced by
$$ 
\sum_{f'(b) \leq n \leq f'(a)}c(n)W(n),
$$ 
where, instead of \eqref{shukou},  $W(n)$ is given by
\begin{equation} \label{shukoubetukei}
W(n)=\frac{1-i}{\sqrt{2}}\frac{\varphi(x_n)}{\sqrt{|f''(x_n)|}}\exp(2\pi i(f(x_n)-nx_n)).  
\end{equation}
\end{remark}

%***********************************************

\section{Proof of Theorem 1}

%**********************************************
Instead of \eqref{approx} we shall use the expression of $Z(t)$ containing a smooth weight $\rho(t)$.
This is because that the error term in \eqref{approx} is too big for our purpose. If we use \eqref{AFE} of 
Lemma \ref{lem-1} we get
\begin{align}  \label{weight-approx}
Z(g(n))=2(-1)^n\sum_{m=1}^{\infty} m^{-1/2}\rho\left(m\sqrt{\frac{2\pi}{g(n)}}\right)\cos(g(n)\log m)+O(g(n)^{-5/6}).
\end{align}

As for the sum of $Z(g(n))$ we consider the sum over even $n$ and odd $n$ separately. First we consider the even $n$ case.
From \eqref{weight-approx} it follows that
\begin{align*}
\sum_{0 \leq n \leq N}Z(g(2n))&=2\sum_{m=1}^{\infty}\frac{1}{m^{1/2}}\sum_{0 \leq n \leq N}
\rho\left(m\sqrt{\frac{2\pi}{g(2n)}}\right)\cos(g(2n)\log m) \\
& \quad +O(N^{1/6}\log^{5/6}N).
\end{align*}
The sum over $m$ is actually a finite sum, in fact $m$ runs over from 1 to $2\sqrt{g(2N)/2\pi}$
and for such $m$, $n$ runs over under the condition $\frac{m^2}{4} \leq \frac{g(2n)}{2\pi} \leq \frac{g(2N)}{2\pi}$.
The contribution from $m=1$ becomes
\begin{equation*}
2\sum_{0 \leq n \leq N} \rho\left(\sqrt{\frac{2\pi}{g(2n)}}\right)=2\left(N+\rho\left(\sqrt{\frac{2\pi}{g(0)}}\right)\right)
\end{equation*}
(note that $\sqrt{2\pi/g(0)}=0.594 \ldots $).

To consider the sum from $m \geq 2$, let
$$
S_m=\sum_{\substack{0 \leq n \leq N \\ \frac{m^2}{4} \leq \frac{g(2n)}{2\pi} }} 
    \rho\left(m\sqrt{\frac{2\pi}{g(2n)}}\right) e^{2\pi i f_m(n)},
$$
where
\begin{equation}  \label{def-f}
f_m(x)=\frac{1}{2\pi}g(2x)\log m.
\end{equation}
Thus we get
$$
\sum_{0 \leq n \leq N}Z(g(2n))=2N+\sum_{m=2}^{\infty}\frac{1}{m^{1/2}}\left(S_m+\bar{S}_m \right) +O(1).
$$
From \eqref{gramfunction} the first and second derivatives of $f_m(x)$ are given by 
\begin{align}
f_m'(x)&=\frac{\log m}{\theta'(g(2x))} \label{f-prime} \\
\intertext{and}
f_m''(x)&=-\frac{2\pi (\log m) \, \theta''(g(2x))}{\theta'(g(2x))^3}, \label{f-doubleprime}
\end{align}
respectively. 
Since $\theta'(t)$ is positive and increasing for $t \geq 6.5$ by \eqref{theta1-seikaku} and \eqref{theta2-seikaku}, $f_m'(x)$ 
is positive and decreasing for $x$ such that $g(2x) \geq  6.5$.

Let $M_0=M_0(m)=m^2/4$ and $M_j=32^jM_0$.  For $M \geq M_0$ we put
\begin{equation*}  %\label{SmM}
S_m(M)=\sum_{M \leq \frac{g(2n)}{2\pi} < 32M}\rho\left(m\sqrt{\frac{2\pi}{g(2n)}}\right) e^{2\pi i f_m(n)}.
\end{equation*}
Then we get the decomposition
$$
S_m=S_m(M_0)+S_m(M_1)+ S_m(M_2) + \cdots + S_m(M_J'),
$$
where the last sum is taken over the range $M_J \leq \frac{g(2n)}{2\pi} \leq \frac{g(2N)}{2\pi}$, $J$ is the 
largest integer such that $32^JM_0 \leq g(2N)/2\pi$ and in fact $J=O(\log N)$.

Now we consider the sum $S_m(M_j)$ in more details. For $m=2$ and $M_0=1$, $S_2(1)$ consists of finite Gram
points, hence $S_2(1)=O(1)$.

To treat other cases, let $h(y)$ be the inverse
function of $g(2x)/2\pi=y$, namely, $h(y)=x$. We have
\begin{align} \label{h-function}
\frac{g(2h(y))}{2\pi}=y \quad \text{and}\quad h\left(\frac{g(2x)}{2\pi}\right)=x.
\end{align}
So the summation condition of $S_m(M)$ is converted to
$$
h(M) \leq n < h(32M).
$$
We also note that $h(y) \sim \frac{y}{2}\log y$.
For the cases other than $m=2$ and $M_0=1$, $f_m'(x)$ is positive and decreasing, and 
from \eqref{f-prime},  \eqref{h-function} and \eqref{theta1-seikaku} we have
\begin{align}  \label{fprime-1}
f_m'(h(M))&=\frac{\log m}{\theta'(g(2h(M)))}=\frac{\log m}{\theta'(2\pi M)} \\
        &=\frac{\log m}{\frac12\log M-\frac{1}{192\pi^2 M^2}+V_1(2\pi M)}.  \nonumber
\end{align}

First we consider the sum $S_m(M_j)$ for $j \geq 1$. For $x \in [h(M_j), h(32M_j)]$, we see
that 
$$
0 <f_m'(h(32M_j)) \leq f_m'(x) \leq f_m'(h(M_j)).
$$
From \eqref{fprime-1} and \eqref{theta1-seikaku} we find that
\begin{align} \label{fprime-3}
f_m'(h(M_j)) &=\frac{\log m}{\frac12\log(8M_{j-1})+\log 2-\frac{1}{192\pi^2 M_{j}^2}+V_1(2\pi M_j)}\\
             &< \frac{\log m}{\frac12\log(8M_{j-1})} < 1  \nonumber
\end{align}
and 
\begin{align}  \label{sitakara}
f_m'(h(32M_j))&=\frac{\log m}{\frac12\log(32M_j)-\frac{1}{192\pi^2 (32M_j)^2}+V_1(64\pi M_j)} \\
              & \geq \frac{2\log m}{\log 32M_j}.  \nonumber 
\end{align}
Here we have used the inequalities
$$
\log 2-\left(\frac{1}{192\pi^2 M_j^2}+|V_1(2\pi M_j)|\right)>0
$$
and
$$
\frac{1}{192\pi^2 (32M_j)^2}-|V_1(64\pi M_j)|>0
$$
for $m \geq 2$ and  $j \geq 1$.  
Since the conditions of Lemma \ref{lem-fdt} are satisfied we can apply \eqref{exp-sum-2} to the sum $S_m(M_j)$ with 
$\delta=\frac{2\log m}{\log 8M_{j-1}}$ and $\delta_1=\frac{2\log m}{\log 32M_j}$.
As a result we get
\begin{equation*}
S_m(M_j) \ll \frac{\log 32M_j}{2\log m}+\frac{\log (M_j/4)}{(5j-4)\log 2} \ll \frac{j}{\log m}+\frac{\log m}{j}+1
\end{equation*}
for $j \geq 1$. Hence putting $l_N=[2\sqrt{g(2N)/2\pi}]$, the contribution of these terms to the sum
$\sum_{2 \leq m \leq l_N} m^{-1/2} S_m$ becomes
\begin{align}  \label{sono1}
\sum_{2 \leq m \leq l_N} \frac{1}{\sqrt{m}}\sum_{1 \leq j \leq J} S_m(M_j) 
&\ll \sum_{2 \leq m \leq l_N} \frac{1}{\sqrt{m}}\left(\frac{\log^2N}{\log m}+\log N+\log m \cdot \log \log N\right)\\[1ex]
& \ll N^{1/4} \log^{3/4}N \log \log N.  \nonumber 
\end{align}

%**********************************************************************************

For $S_m(M_0)$ $(m \geq 3)$, we apply Lemma \ref{lem-3} with $\varphi(x)= \rho(m \sqrt{2\pi/g(2x)})$
and $f(x)=f_m(x)$ for $x \in [h(32M_0), h(M_0)]$.  For this we have to check the assumptions in Lemma \ref{lem-3}.
By \eqref{f-doubleprime}, \eqref{theta1-asymp} and \eqref{theta2-asymp} we have
%\begin{align*}
%f_m''(x)=-\frac{2\pi \log m \, \theta''(g(2x))}{\theta'(g(2x))^3},
%\end{align*}
%hence
$|f_m''(x)| \asymp \frac{\log m}{M_0 \log^3 M_0}$. Furthermore we have easily that
$$
f_m'''(x) \ll \frac{\log m}{M_0^2 \log^4M_0}, \qquad f_m^{(4)}(x) \ll \frac{\log m}{M_0^3 \log^5 M_0}.
$$
Hence we can take $A=\frac{M_0 \log^3 M_0}{\log m}, \quad U=M_0 \log M_0$. 
As we remarked in Remark \ref{remark-KV}, the assertion of Lemma \ref{lem-3} is valid though the condition $A \ll U$ is not satisfied.
On the other hand, we have
\begin{align*}
\varphi'(x)=-\rho'\left(m\sqrt{\frac{2\pi}{g(2x)}}\right)\sqrt{2\pi}m g(2x)^{-3/2}g'(2x).
\end{align*}
Since $g(2x) \asymp M_0$, we have $\varphi'(x) \ll \frac{1}{M_0 \log M_0}= U^{-1}.$  Similarly we see that
$\varphi''(x) \ll U^{-2}$ holds. (Here we can take $H = c \max\{\rho(x), |\rho'(x)|, |\rho''(x)|\} \ll 1$, where $c$ is a constant.)
Furthermore we have to note that $f_m''(x)$ is negative in our case,
hence the sum on the right hand side of \eqref{exp-sum-trans} should be taken in the range 
$f_m'(h(32M_0)) \leq \nu \leq f_m'(h(M_0))$.
Now by \eqref{fprime-1} we have
\begin{align}  \label{fprime-2}
f_m'(h(M_0))&=\frac{\log m}{\log m-\left(\log 2+\frac{1}{12\pi^2m^4}-V_1(\frac{\pi m^2}{2})\right)}.
\end{align}
Similarly as above, if we use \eqref{theta1-seikaku} in Lemma \ref{lem-4} we find that 
\begin{align*}
2<f_m'(h(M_0)) <3 & \quad \text{for $m=3, 4$} \\
\intertext{and}
1 <f_m'(h(M_0)) < 2 & \quad \text{for all $m \geq 5$.}
\end{align*}
On the other hand we have already seen in \eqref{fprime-3} that
\begin{align*} %\label{fprime-32M0} 
f'(h(32M_0)) <\frac{\log m}{\log m+\frac12 \log 2}<1.
\end{align*}
Thus we get
\begin{align*}
S_m(M_0)&=W_m(1)+R_m   \quad \mbox{for $m \geq 5$,} \\
S_m(M_0)&=W_m(1)+W_m(2)+R_m \quad \mbox{for $m =3,4$}, 
\end{align*}
%The above consideration shows that 
%\begin{enumerate}
%\item[(i)]  for any $m \geq 5$, the interval $[f_m'(h(32M_0)), f_m'(h(M_0))]$ contains only one integer $\nu=1$.
%\item[(ii)] for $m=3$ and $4$, the interval $[f_m'(h(32M_0)), f_m'(h(M_0))]$ contains only two integers $\nu=1$ and $2$.
%\end{enumerate}
where 
\begin{equation}  \label{W-m}
W_{m}(\nu)=\frac{1-i}{\sqrt{2}}\cdot \frac{\rho\left(m\sqrt{2\pi/g(2x_\nu)}\right)}{\sqrt{|f_m''(x_\nu)|}}
\exp\Bigl(2\pi i(f_m(x_\nu)-\nu x_\nu)\Bigr)
\end{equation}
(\eqref{shukoubetukei} must be used since $f''_m(x)$ is negative) and
\begin{equation*}  %\label{gosa}
R_{m} \ll \log m+T_{h(32M_0)}+T_{h(M_0)}+\log(f_m'(h(M_0))-f_m'(h(32M_0))+2)
\end{equation*}
(see \eqref{gosaR} and \eqref{def-T} for the definitions of $R_m$ and $T_{\mu}$).
The contribution of $S_m(M_0)$ to the sum $\sum_{2 \leq m \leq l_N} m^{-1/2}S_m$ is 
\begin{equation}  \label{saishuu}
\sum_{3 \leq m \leq l_N}\frac{W_m(1)}{\sqrt{m}}+\sum_{3 \leq m \leq l_N}\frac{R_m}{\sqrt{m}}+O(1),
\end{equation}
which we shall calculate now.

First we shall determine the explicit form of $W_{m}(1)$ given by \eqref{W-m}.
Define $x_1$ so as to satisfy $f_m'(x_1)=1$, which is equivalent to $\theta'(g(2x_1))=\log m$.
Using \eqref{theta1-asymp} we obtain
\begin{equation} \label{shiki-1}
\frac{g(2x_1)}{2\pi}\left(1+O\left(\frac{1}{g(2x_1)^2}\right)\right)=m^2.
\end{equation}
On the other hand, by the definition of $g(2x)$ we have $\theta(g(2x_1))=2\pi x_1$, hence
$$
\theta\left(2\pi m^2+ O\left(\frac{1}{g(2x_1)}\right)\right)=2\pi x_1.
$$
Using \eqref{theta} on the left hand side of the above formula, we obtain
\begin{align}  \label{shiki-2}
x_1=m^2\log m - \frac12m^2-\frac{1}{16}+O\left(\frac{\log m}{m^2}\right).
\end{align}
Therefore from \eqref{shiki-1} and \eqref{shiki-2} we get
\begin{align*} % \label{shiki-3}
f_m(x_1)-x_1&=\frac{g(2x_1)}{2\pi} \log m -x_1 \\
          &=\frac12 m^2+\frac{1}{16}+O\left(\frac{\log m}{m^2}\right).  \nonumber
\end{align*}
Furthermore from \eqref{theta1-asymp}, \eqref{theta2-asymp} and \eqref{shiki-1} we have
\begin{align*} %\label{shiki-4}
\theta'(g(2x_1))=\log m+O\left(\frac{1}{m^4}\right) 
\intertext{and}
\theta''(g(2x_1))=\frac{1}{4\pi m^2}+O\left(\frac{1}{m^6}\right).
\end{align*}
Hence
\begin{align*} %\label{shiki-5}
f_m''(x_1)&=-\frac{2\pi \log m \, \theta''(g(2x_1))}{\theta'(g(2x_1))^3}
=-\frac12 \frac{1}{(m \log m)^2}\left(1+O\left(\frac{1}{m^{4}}\right)\right).
\end{align*}
Combining these, we get
\begin{align*}
W_{m}(1)&=e^{-\frac{\pi}{4}}\rho\left(m\sqrt{\frac{2\pi}{g(2x_1)}}\right)\sqrt{2}m \log m \left(1
+O\left(\frac{1}{m^{4}}\right)\right) \\
      & \qquad \times e^{2\pi i \left(\frac{m^2}{2}+\frac{1}{16}+O\left(\frac{\log m}{m^2}\right)\right)}  \\
      &=\sqrt{2}e^{-\frac{\pi i}{8}}\rho\left(m\sqrt{\frac{2\pi}{g(2x_1)}}\right) (-1)^m m \log m
         \left(1+O\left(\frac{\log m}{m^2}\right)\right).
\end{align*}
Now from \eqref{shiki-1} and \eqref{rhokyodo} of Lemma \ref{lem-2}, we have
\begin{align*}
\rho\left(m\sqrt{\frac{2\pi}{g(2x_1)}}\right)=\rho(1+O(m^{-2}))=\frac12+O\left(m^{-C} \right)
\end{align*}
for any large $C>0$. Hence we get
\begin{align*}  %\label{shukou-1}
W_{m}(1)=\frac{1}{\sqrt{2}}e^{-\frac{\pi i}{8}}(-1)^m m \log m + O\left(\frac{\log m}{m}\right).
\end{align*}
%Similarly to the above case we have
%\begin{align}
%W_{m,M_0}(2)=e^{-\frac{\pi}{4}}\rho(...) \frac{\sqrt{m}\log m}{2\sqrt{2}} e^{2\pi i (m+\frac{1}{16})}\left(1+O\left(\frac{\log m}{m}\right)\right).
%\end{align}
%But as is seen by Lemma 5, these terms only appear for m=3 and 4, hence their contribution becomes finite.
Therefore we find that 
%Now we can calculate $Q_1$ in \eqref{lastform}. Using \eqref{shukou-1} we see that 
\begin{align}  \label{sono2}
\sum_{3 \leq m \leq l_N} \frac{W_m(1)}{\sqrt{m}}
&=\frac{1}{\sqrt{2}}e^{-\frac{\pi i}{8}}\sum_{3 \leq m \leq l_N}(-1)^m m^{1/2} \log m +O(1) \\
&=\frac{1}{2\sqrt{2}}e^{-\frac{\pi i}{8}}(-1)^{l_N}l_N^{1/2}\log l_N +O(1) \nonumber \\[1ex]
& \ll g(2N)^{1/4} \log g(2N) \nonumber \\[1ex]
& \ll N^{1/4} \log^{3/4} N. \nonumber
\end{align}
The second equality is obtained as follows.
Suppose first that $l_N=2L$ (even integer). Then we have
\begin{align*}
U:&=\sum_{1 \leq m \leq 2L}(-1)^m m^{1/2}\log m  \\
&=\sum_{k=1}^L\left((2k)^{1/2}\log 2k- (2k-1)^{1/2}\log (2k-1)\right).
\end{align*}
Since
\begin{align*}
&(2k)^{1/2}\log 2k- (2k-1)^{1/2}\log (2k-1) \\
& \qquad =\left(\frac{1}{\sqrt{2}}+\frac{\log 2}{2\sqrt{2}}\right)\frac{1}{k^{1/2}}+\frac{1}{2\sqrt{2}}\frac{\log k}{k^{1/2}}
+O\left(\frac{\log k}{k^{3/2}}\right), \\
& \sum_{1 \leq k \leq L} \frac{1}{k^{1/2}}=2L^{\frac12}+O(1) \\
\intertext{and} 
&\sum_{1 \leq k \leq L} \frac{\log k}{k^{1/2}}=2L^{\frac12}\log L -4 L^{\frac12}+O(1),
\end{align*}
we get
$$
U=\frac12(2L)^{\frac12}\log (2L) +O(1).
$$
The case $l_N=2L+1$ (odd integer) is similar.

Next we treat the contribution from $R_m$ in \eqref{saishuu}. It is enough to consider the sum for $m \geq 5$. 
For such $m$, we know that $f'(h(32M_0)) < 1 < f'(h(M_0)) <2$ and from \eqref{sitakara} (which is also true for $j=0$) 
and \eqref{fprime-2} we have $T_{f'(h(32M_0))}, T_{f'(h(M_0))}  \ll \log m$ and hence
$$
R_m \ll \log m.
$$
Therefore we find  that 
\begin{equation} \label{sono3}
\sum_{3 \leq m \leq l_N} \frac{R_m}{\sqrt{m}} \ll N^{1/4} \log^{3/4} N.
\end{equation}
From \eqref{sono1}, \eqref{sono2} and \eqref{sono3} we get the assertion \eqref{evencase}.

%
%Now we can calculate the sum of $S_m$. Putting  $l_N=2\sqrt{g(2N)/2\pi}$ we have
%\begin{align*}
%\sum_{2 \leq m \leq l_N} \frac{S_m}{\sqrt{m}}
%&= \sum_{3 \leq m \leq l_N}\frac{W_{m,M_0}(1)}{\sqrt{m}}
% + \sum_{3 \leq m \leq l_N} \frac{1}{\sqrt{m}}\sum_{M \geq M_0}R_{m,M} \\
%&  \quad + \frac{1}{\sqrt{2}}\sum_{M \geq 8}R_{2,M} + O(1) \\
%%    %W_3(2)+W_4(2)+S_2(1).
%&= Q_1+Q_2+Q_3+O(1),
%\end{align*}
%say, where the symbol $\sum_{M \geq M_0}$ means that $M$ runs over $M_0, 8M_0, 64M_0, \ldots$.
%From \eqref{shukou-1}, we have
%\begin{align}
%Q_1&=\frac{1}{\sqrt{2}}e^{-\frac{\pi i}{8}}\sum_{3 \leq m \leq l_N}(-1)^m m^{1/2} \log m +O(1) \\
%&=\frac{1}{\sqrt{2}}e^{-\frac{\pi i}{8}}(-1)^{[l_N]}[l_N]^{1/2}\log [l_N] +O(1) \nonumber \\
%& \ll g(2N)^{1/4} \log g(2N) \ll N^{1/4} \log^{3/4} N,  \nonumber
%\end{align}
%where $[l_N]$ is a greatest integer which does not exceed $l_N$.
%
%Similarly from \eqref{gosakou-1} and \eqref{gosakou-2}, we have
%$$
%Q_2 \ll N^{1/4} \log^{7/4}N, \quad Q_3 \ll \log^3 n.
%$$
%This proves the assertion \eqref{evencase}.

The assertion \eqref{oddcase} is proved similarly.

\qed

\bigskip

\begin{flushleft}
Xiaodong Cao \\
Department of Mathematics and Physics, \\
Beijing Institute of Petro-Chemical Technology,\\
Beijing, 102617, P. R. China \\
e-mail: caoxiaodong@bipt.edu.cn

\medskip

Yoshio Tanigawa\\
Nishisato 2-13-1, Meitou, Nagoya 465-0084, Japan\\
e-mail:  tanigawa@math.nagoya-u.ac.jp,   tanigawa\verb'_'yoshio@yahoo.co.jp 

\medskip

Wenguang Zhai\\
Department of Mathematics, \\
China University of Mining and Technology, \\
Beijing 100083, P. R. China\\
e-mail: zhaiwg@hotmail.com
\end{flushleft}

\end{document}